\documentclass[11pt, a4paper]{amsart}
\usepackage{amsfonts}
\usepackage{amsmath}
\usepackage{epsfig}
\usepackage{graphicx}
\usepackage{psfrag}
\newcommand{\equ}[1]{(\ref{#1})} 
 
\newcommand{\be}{\begin{equation}} 
\newcommand{\ee}{\end{equation}} 
\newtheorem{lem}{Lemma}[section] 
\newcommand{\rr}{{\mathbb R}} 
\newcommand{\nn}{{\mathbb N}} 
\newtheorem{teo}{Theorem}[section]

\newtheorem{rem}{Remark}[section] 
\newtheorem{cor}{Corollary}[section] 
\newcommand{\calA}{{\mathcal A}}

\newcommand{\calH}{{\mathcal H}}

\newcommand{\calT}{{\mathcal T}}

\begin{document} 
\bibliographystyle{plain} 
\title[Relaxation of the flow of triods]{ Relaxation of the flow of triods by Curve Shortening Flow via the
  vector-valued parabolic Allen-Cahn equation}

\author{Mariel S\'aez Trumper}
\address{Mariel S\'aez Trumper
\hfill\break\indent
Max Planck Institute for Gravitational Physics
\hfill\break\indent
Albert Einstein Institute
\hfill\break\indent
Am M{\"u}hlenberg 1\\
D-14476 Golm\\
\hfill\break\indent
Germany .}
\email{{\tt  mariel.saez@aei.mpg.de}}

\begin{abstract}
In this paper we find solutions $u_\epsilon$ to a certain class of  vector-valued parabolic Allen-Cahn equation that as $\epsilon \to 0$ develops as interface a given triod evolving under curve shortening flow. 
\end{abstract}

\maketitle

\section{Introduction}

This paper studies the relationship between a vector valued Allen-Cahn equation and the motion of triods by mean curvature flow. Here we exploit the techniques previously developed in \cite{tesis}, \cite{papertesis}.

The Allen-Cahn Equation is given in a domain $\Omega$ by:
\begin{align}\frac{\partial u_\epsilon}{\partial t}-\Delta u_\epsilon+ \frac{\nabla_uW(u_\epsilon)}{\epsilon^2}&=0 \hbox{ for } x\in \Omega \label{laeq}\\
 u_\epsilon(x,0)&=\psi_\epsilon(x),\label{ci} \\
u_\epsilon|_{\partial \Omega}&=\phi_\epsilon (x,t) \label{bc} \end{align}
where $u_\epsilon:\rr^n\times \rr_+\to \rr^m$ and $W:\rr^n\to \rr$ is a positive potential with a finite number of minima. In particular we will concentrate on the case $m=n=2$ and $W$ a function with 3 minima. We prove that triods evolving under curve shortening flow can be realized as nodal sets of this equation
(for a precise statements and definitions see Section \ref{not}). We also include
 some corollaries derived from this representation (which are  stated in section  \ref{not}, as well).

Equation \equ{laeq} has been studied by several authors. In particular the scalar case is widely known (that is when $m=1$). It has been shown 
for  double well potentials that scalar solutions  $u_\epsilon$ to \equ{laeq} converge as $\epsilon\to 0$ almost everywhere to minima of $W$ and they  develop interfaces separating the regions where the minima are attained. These interfaces  evolve under mean curvature flow. See  
\cite{froprogb, quaoptgb},  \cite{genandxc, geoevpd,phatrale}, \cite{papertesis} for precise statements.

In the vector-valued case less is known. Some results for the stationary equation can be found in 
\cite{stalasa,  anaofna, onthena,expstana},
 \cite{minintsb}, \cite{atrhlb}, \cite{vecvalps} and \cite{stationary}. For the parabolic problem  L.Bronsard and F.Reitich(\cite{onthrlb}) predicted, via a formal analysis, that for a 3 well potential solutions $u_\epsilon$ to \equ{laeq} converge almost everywhere to minima of the function $W$ and that the develop interfaces evolve under curve shortening 
flow. 
In particular, Bronsard and Reitich~\cite{onthrlb} conjectured that as $\epsilon\to 0$ the solutions might develop  a {\em triod} structure. That is the nodal set (or interface set)  is a network composed by three regular curves which meet at a unique point and each of them evolves under curve shortening flow. In this paper we give a rigorous proof of this fact for a certain class of potentials $W$.

The flow of triods under curve shortening flow can be described analytically by equation \equ{triod} in  Section \ref{not}. In order to have a well defined system of equations, an extra condition is necessary at the meeting point. In \cite{onthrlb} arbitrary prescribed angles were considered.
 In this context the authors
proved short-time existence of triods, under the assumption that the initial condition satisfies strong compatibility conditions. 
Recently great progress was made by Mantegazza, Novaga and Torterelli~\cite{moybycm} for triods that meet at $120^0$ angles for every $t\geq 0$.  They  proved long time existence (up to the first singularity time) 
with generic initial data  satisfying the meeting condition at the triple point.
In particular they were able to remove the compatibility conditions, but they were not able to prove geometric uniqueness in the more general case. In this paper, using the representation provided by equation \equ{laeq}, we show that indeed geometric uniqueness holds (see Corollary \ref{cor1}).

Recently, O.Schn\"urer and F. Schulze \cite{selsimos} considered 
 triods evolving under curvature flow and meet at $120^0$ for every $t> 0$, 
but do not necessarily satisfy this condition at $t=0$. They showed that when  3 lines meeting at any arbitrary angles are considered as initial condition, there is a self-similar solution to the Curve Shortening Flow equation for triods ( see equation \equ{triod} in Section \ref{not}), such that for every $t>0$ the three curves  meet at $120^0$ angles.
 It is expected in the general case (i.e. any initial condition meeting at arbitrary angles is considered) that these self-similar solutions will predict the behavior of the triple point for short-time. As first step, in 
Corollary \ref{cor2} we can show in certain situations that  the backward blow up at the triple point of  smooth solutions to the flow correspond to one of the self-similar solutions described in \cite{selsimos}.

%The results of this paper only consider times before the first singularity occ%urs.

We organize this paper as follows: in section \ref{not} we establish some notation and we make precise the statements of the Theorem and its Corollaries. We also include in that section the statement of some lemmas that we will use but were already  proved in the literature. In sections \ref{pfth} and \ref{pfcor} we include the proofs of the main Theorem and its Corollaries respectively.
We finish with Section \ref{finsec}, where  we state some  problems that remained open.

\section{Notation and Results} \label{not}

Consider $\Omega\subset \rr^2$ an open domain in $\rr^2$ with smooth boundary.
We say that a {\em triod} $\calT=\{\gamma^i:[0,1]\times[0,T)\to 
\Omega\}_{i=1}^3$ of curves evolves under curve shortening flow (with Dirichlet Boundary data)
if it satisfies
the following system of equations:
\be\left\{\begin{array}{ll}
\gamma^i_\lambda(\lambda,t)\ne 0& \lambda \in [0,1], t \in [0,T], i\in\{1,2,3\}  \\
\gamma^i(\lambda_1,t)\ne \gamma^i(\lambda_2,t) \hbox{ if } \lambda_1\ne \lambda_2&  \lambda_1, \lambda_2 \in [0,1], t \in [0,T], i\in\{1,2,3\}\\
\gamma^i(\lambda_1,t)=\gamma^j(\lambda_2,t)& \hbox{ iff } \lambda_1=\lambda_2=0  \hbox{ and }i\ne j \in \{1,2,3\}
\\
%\sum_{i=1}^3 \frac{\gamma_\lambda^i(0,t)}{|\gamma_\lambda^i(0,t)|}=0 \hbox{ for }t>0& \\
\gamma^i(1,t)=P^i(t)\in \partial \Omega &  \\
\gamma^i(\lambda,t)=\sigma^i(\lambda) & \\
\gamma^i_t(\lambda,t)=\frac{\gamma^i_{\lambda\lambda}(\lambda,t)}{|\gamma^i_\lambda(\lambda,t)|^2}& \hbox{ in } \Omega.  
\end{array}\right. \label{triod} \ee

Suppose that there is a solution to \equ{triod}  $\{\gamma^i\}_{i=1}^3$ satisfying  an appropriate condition  for the  angles between the curves at the triple point (we are going to consider fixed angles for all times). Moreover, assume there is a potential $W:\rr^2\to \rr$ that is consistent with this angle condition (this consistency of $W$ as well as the other necessary conditions on the potential are  going to be specified later in this section).
In this paper we show that exists an  initial condition $\psi_\epsilon$, a boundary condition $\phi_\epsilon$ and a solution $u_\epsilon$
to \equ{laeq}-\equ{ci}-\equ{bc}
which nodal set as $\epsilon\to 0$ agrees with $\{\gamma^i\}_{i=1}^3$
In particular, we need to require that $W:\rr^2\to \rr$ 
is a $C^3$ function that satisfies
\begin{enumerate}
\item[(W1)] \label{condW0} $W$ has only three local minima $c_1, c_2$ and $c_3$ and $W(c_i)=0$; 
\item[(W2)] \label{cond1W} the matrix $\frac{\partial^2 W(u)}{\partial u_i \partial u_j}$ is positive definite at $\{c_i\}_{i=1}^3$, that is the minima are nondegerate;
\item[(W3)] \label{cond2W} there exist positive constants $K_1, K_2$ and $m$, and a number $p\geq 2$ such that
$$K_1|u|^p\leq W(u)\leq K_2 |u|^p \hbox{ for } |u|\geq m;$$
\item[(W4)] \label{cond30W} $V(r,\theta):=W(u+r(\cos \theta, \sin \theta))=r^2+O(r^3)$
for $r$ sufficiently small and $u=c_i$ for some $i\in\{1,2,3\}$, where $r$ and $\theta$ are local polar coordinates.

\item[(W4)] \label{cond3W} There is a $K>0$ such that $\frac{\partial^2 W(u)}{\partial u_i \partial u_j}$ is positive definite for every $u>K$.
\end{enumerate}

Results in \cite{vecvalps} and \cite{onthrlb}  suggest that under these hypothesis  
the following stationary solutions 
to \equ{laeq} exist:
\begin{itemize}
\item \be \zeta_{ij}''(\lambda)+\frac{\nabla W(\zeta_{ij}(\lambda))}{2}=0 \label{eczeta}.\ee
\be \lim_{\tau\to -\infty}\zeta_{ij}(\tau)=c_i, 
\lim_{\tau\to \infty}\zeta_{ij}(\tau)=c_j, \label{explim}\ee
where these limits are attained at an exponential rate. A standard argument implies that this convergence rate  also holds for the derivatives of $\zeta_{ij}$.

Recent work of N. Alikakos, S.Betel\'u and X.Chen~\cite{expstana} proves that this is not always the case. In this paper we {\bf assume} the existence of the curves $\zeta_{ij}$ satisfying \equ{eczeta}-\equ{explim}. Under this condition results in \cite{stationary} imply:

\smallskip 

\item There is a  stationary solution to \equ{laeq} $u_*(x):\rr^2\to \rr^2$  that satisfies as $r\to \infty$
\be u_*(r \cos\theta,r\sin\theta) \to c_i \hbox{ for } \theta \in 
[\theta_{i-1}, \theta_i], \label{convu*} \ee
\be  u_*(r \cos\theta_i,r\sin\theta_i)\to \zeta_{ij}(0). \label{convu*2}\ee
where $\theta_i$ are given by the function $W$ in the following manner:

Define
 \begin{align} \Gamma (\zeta_1,\zeta_2)=\inf\left\{\int_0^1W^{\frac{1}{2}}(\gamma(\lambda))|\gamma'(\lambda)|d\lambda: \right.&
 \gamma \in C^1([0,1],\rr^2), \notag \\
&\left. \gamma(0)=\zeta_1\hbox{ and }\frac{}{} \gamma(1)=\zeta_2 \right\}. \label{defdist}\end{align}
Consider $\{\alpha_i\}_{i=1}^3\in[0,2\pi)$ such that
\be \frac{\sin \alpha_1}{\Gamma(c_2,c_3)}= \frac{\sin \alpha_2}{\Gamma(c_1,c_3)}= \frac{\sin \alpha_3}{\Gamma(c_1,c_2)}.\label{condang}\ee
Then the angles $\theta_i\in [0, 2\pi)$ are uniquely determined by $\alpha_i=\theta_{i+1}-\theta_i$.

\begin{rem}\label{dec}We would like to remark that the by \cite{stationary} the  convergences in equations \equ{convu*} and \equ{convu*2}  (as well as convergence of the derivatives of $u_*$) 
are of order $r^{-m}$ for every $m>0$. That is, for every
 $m>0$ and $n\in \nn$ there is a constant $C$ (that might depend on $m,n$) such that
\begin{align} \left|\frac{\partial ^n}{\partial x_1^j \partial x_2^{n-j}} u_*(r \cos\theta,r\sin\theta)- c_i\right|\leq & \frac{C}{1+r^m} \label{convu*dec} \\& \hbox{ for } j\leq n, 
 \quad \theta \in 
(\theta_{i-1}+\delta, \theta_i-\delta), \notag \end{align}

\begin{align}  \left|\frac{\partial ^n}{\partial x_1^j \partial x_2^{n-j}}  u_*(r \cos\theta,r\sin\theta)- \zeta_{ij}\left(d_i(x)\right)\right|\leq&  \frac{C}{1+r^m} \label{convu*dec2} \\ &\hbox{ for } j\leq n \hbox{ and }
 \theta \in 
(\theta_{i}-\delta, \theta_i+\delta) \notag \end{align}
where $\delta>0$ is a small enough constant and $d_i$ is the distance to the line of slope $\tan \theta_i$. 
\end{rem}

We  say that the potential $W$ is ``symmetric'' when 
$\Gamma(c_i,c_j)=\Gamma(c_j,c_k)$ for every $i\ne j$, $j\ne k$. Notice that this kind of potential correspond to equal angles $\theta_i=120^0$.
\end{itemize}

We also need some notation and basic definitions for triods that we list in what follows:
\begin{itemize}
\item Let $O(t):=\gamma^i(0,t)$.
\item The sub-(and super)indexes will always be considered modulo 3 in this paper.
\item We denote by $\tau^i(t)=(\cos \theta^i(t),\sin \theta^i(t))$ the tangent to the curve $\gamma^i(x,t)$ at $x=0$ (that is at the meeting point $O(t)$).
In general, we will require $\theta^i(t)$ to be determined by 
$\theta^1(t)$, hence, for simplicity we denote simply $\theta(t)\equiv \theta^1(t)$.
\item Consider  $\tilde{\delta}>0$. Then we say that a triod is {\bf  graphical over $\{\tau^i\}_{i=1}^3 $ } inside the ball $B_{\tilde{\delta}}(O(t))=
\{x\in \rr^2:|x-O(t)|\leq \tilde{\delta}\}$ if the curve $\gamma^i$ can be written as a graph over the line with slope $\tau^i$ that passes through $O(t)$.
 
\item Let $S_{ij}(t)\subset \Omega$ be the region bounded by $\gamma^i(\cdot,t)$ and $\gamma^j(\cdot,t)$.

\item Let $d_i(x,t)=dist(x,\gamma^i(\cdot,t))$ the signed distance of  a point $x\in \rr^2$ to the curve $\gamma^i(\cdot,t)$. 

Since we  consider signed distances, we need to choose a consistent convention regarding the signs. The reader should keep in mind figure 1 for the choice of signs described below.

 \begin{center}
 \psfrag{g1}{$\gamma^1$}
 \psfrag{g2}{$\gamma^2$}
 \psfrag{g3}{$\gamma^3$}
 \psfrag{s12}{$S_{12}$}
 \psfrag{s23}{$S_{23}$}
 \psfrag{s31}{$S_{31}$}
 \psfrag{om}{$\Omega$ }
 \psfrag{fig.1}{Figure 1 }
 \epsfysize=5cm
 \epsfbox{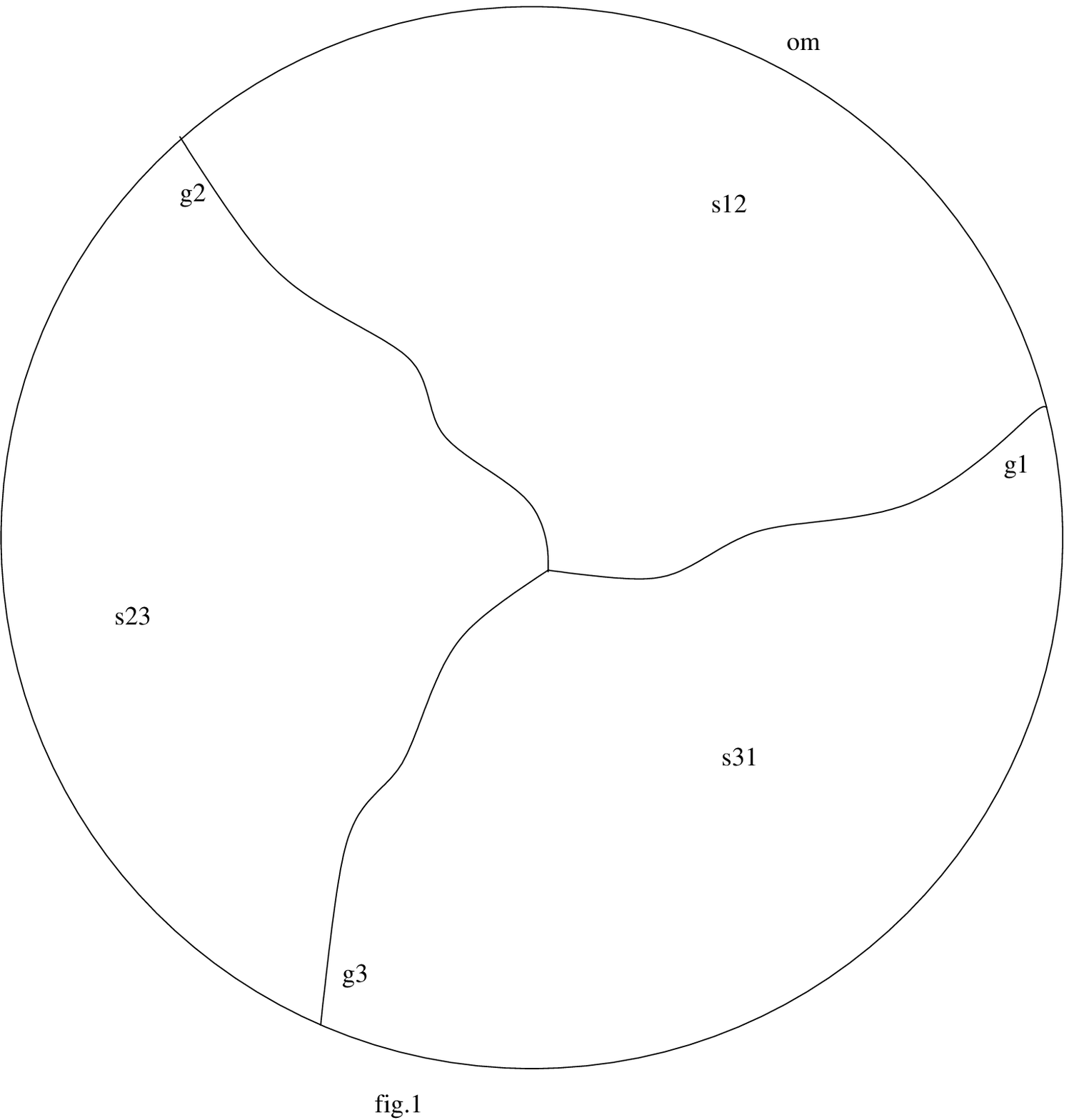}
 \end{center}

We choose $d_i(\cdot,t)$ such that in each sector $S_{i i+1}(t)$ holds $d_i(x,t)>0$ and $d_{i+1}(x,t)<0$ (e.g. $d_1(x,t)>0$ and $d_2(x,t)<0$ for $x\in S_{12}(t)$.)

\item Define $d_{ij}(x,t)$ to be the signed distance to the curve defined 
by $\gamma^i(x,t) \bigcup \gamma^j(x,t)$ where the sign of the distance is chosen to be consistent with $d_i(x,t)$ (e.g. $d_{13}(x,t)>0$ for $x\in S_{12}(t)\bigcup S_{23}(t)$ and $d_{31}(x,t)=-d_{13}(x,t)$.)

\begin{rem}\label{nbang}

Notice that by making $\tilde{\delta}$ small enough we can assure that 
$\gamma^i(\cdot,t)\bigcap B_{\tilde{\delta}}(O(t))=\{\gamma^i(\lambda, t):\lambda\in[0,\lambda_i]\}$ for some $\lambda_i>0$. We want to avoid that curves ``re-enter'' the ball  $B_{\tilde{\delta}}$. For example, in the following picture: 
\newpage

% \begin{figure}
 \begin{center}
 \psfrag{g1}{$\gamma^1$}
 \psfrag{g2}{$\gamma^2$}
 \psfrag{g3}{$\gamma^3$}
 \psfrag{t1}{$\tau^1$}
 \psfrag{t2}{$\tau^2$}
 \psfrag{t3}{$\tau^3$}
 \psfrag{Bdel1}{$B_{\delta_1}(O(t))$ }
 \psfrag{Bdel2}{$B_{\delta_2}(O(t))$ }
 \psfrag{O(t)}{$O(t)$}
\epsfysize=5cm
  \epsfbox{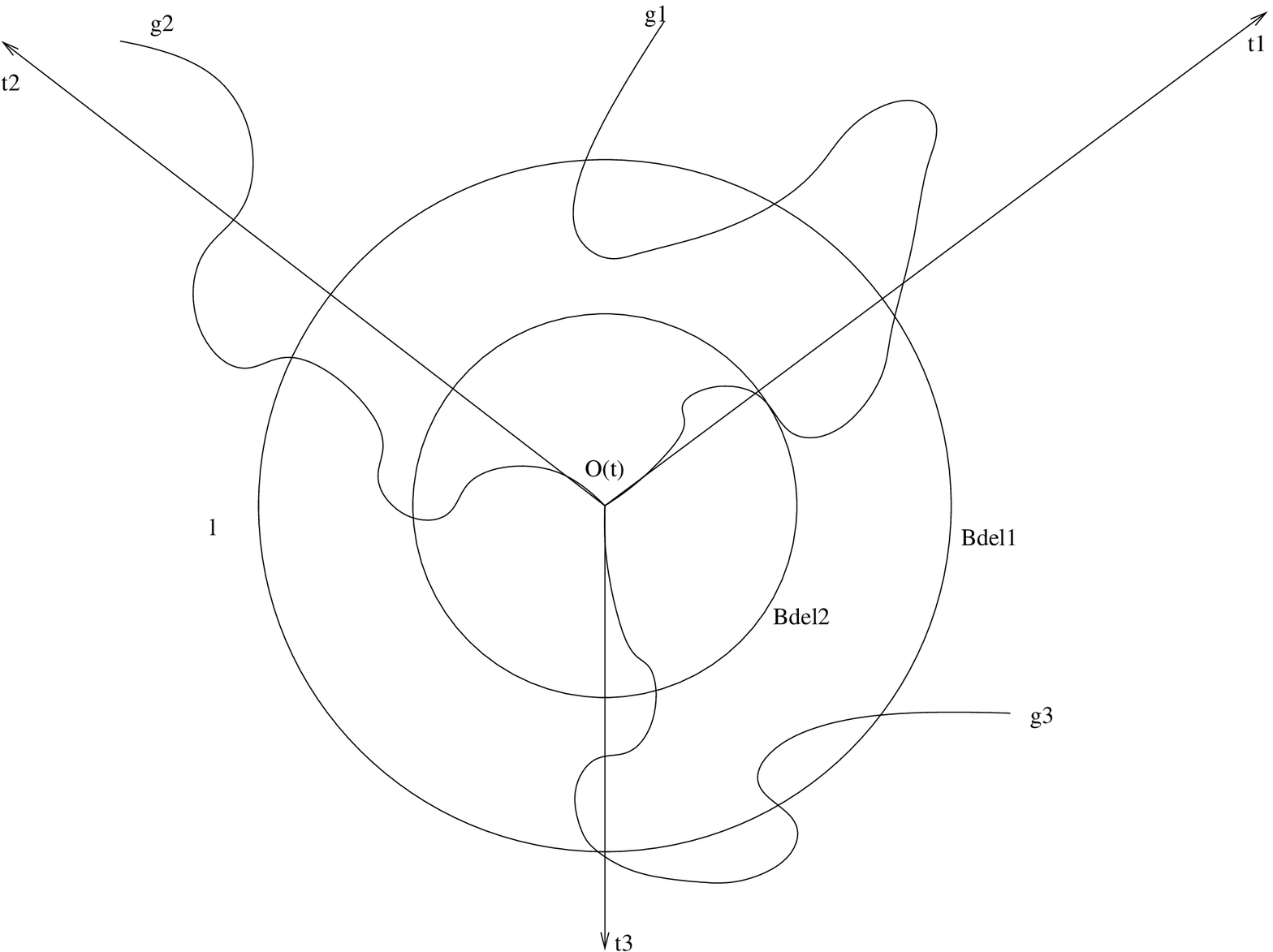}
\end{center}
%\caption{Figure 2}
%\end{figure}

we {\bf will not} want to choose $\tilde{\delta}=\delta_1$, but $\tilde{\delta}=\delta_2.$

Hence, we will assume in what follows that the $\tilde{\delta}$ is always chosen small enough.

In particular, for $\tilde{\delta}$ small enough holds that
$\calT$ is graphical over $\{\tau^i\}_{i=1}^3 $  inside the ball $B_{\tilde{\delta}}(O(t))$. Then, since $\gamma^i(\cdot,t)\bigcap \gamma^j(\cdot,t)=\{O(t)\}$ for $j\ne i$, there is a $\delta_{int}>0$ such that 
\begin{align}\gamma^i(\cdot)\bigcap B_{\tilde{\delta}}(O(t))\subset \left\{O(t)+r(\cos(\theta ^i(t)+\theta),
 \frac{}{} \right.
& \sin(\theta ^i(t)+\theta) )\in \Omega:   \label{gbdeq} \\ & \left.  \theta\in \left(-\frac{\delta_{int}}{2}, \frac{\delta_{int}}{2} \right), r\geq 0 \right\}. \notag \end{align}

\end{rem}

\item Let $R_\theta$ denotes the rotation matrix by an angle $\theta$.

\end{itemize}

Suppose now that we have a triod $\calT$ such that for every $t\in[0,T]$ there is a 
$\tilde{\delta}$ such that is graphical over $\{\tau^i\}_{i=1}^3 $  inside the ball $B_{\tilde{\delta}}(O(t))$ and consider $\delta_{int}$ like in Remark \ref{nbang}.
This allow us to define the following functions:
\begin{itemize}
 \item Let
$\eta_1:\rr^2\to \rr$  such that $\eta_1(x)\equiv 1$ when $|x|\leq\frac{1}{2}$ 
and $\eta_1(x)\equiv 0$ for $|x|\geq 1$.
\item Let
$\eta_2:\rr\to \rr$  such that $\eta_2(x)\equiv 1$ when $|x|\leq\frac{1}{2}$ 
and $\eta_2(x)\equiv 0$ for $|x|\geq 1$
\item Consider a partition of unity $\{\xi^{int}_i\}_{i=1}^6$ associated to  the family of intervals $\{\calA_j\}_{j=1}^6$, where 
%$\{(2\pi-\delta,2\pi]\bigcup[0,\delta), (\frac{\delta}{2},\theta_1-\frac{\delta}{2}),(\theta_1-\delta,\theta_1+\delta),
%(\theta_1+\frac{\delta}{2},\theta_2-\frac{\delta_{int}}{2}),(\theta_2-\delta_{int},\theta_2+\delta_{int}),(\theta_2+\frac{\delta_{int}}{2},2\pi-\frac{\delta_{int}}{2})\}$,
$$\calA_{2i}=\left(\theta_i-\delta_{int}, \theta_i+\delta_{int}\right)$$
$$\calA_{2i+1}=\left(\theta_i+\frac{\delta_{int}}{2},\theta_{i+1}-\frac{\delta_{int}}{2}\right),$$
 %that is $\eta_j(x)\geq 0$, $\eta_j (x)\equiv 0$ for $x\not \in \calA_j$ and
%$\sum_j\eta_j(x)=1$ for every $x$. 
%Since the curves $\gamma_i$ do not intersect each other, it is possible to choose $\delta$ such that $$\gamma_i\in \left\{x=r(\cos(\theta (t)+\theta),\sin(\theta (t)+\theta) )\in \Omega:  \theta\in \left(-\frac{\delta}{2}, \frac{\delta}{2} \right) \right\}.$$
and $\theta_i$ are the angles given by \equ{convu*}.
\item Outside the ball $B_{\tilde{\delta}}(O(t))$ we choose $\delta_i$ such that $d_i(x,t)$ is smooth and well defined for 
$x$ such that $d_i(x,t)\leq \delta_i$ and  satisfies $d_i(x,t)\geq \delta_i$ for every $x\in S_{i+1 i-1}\bigcap (B_{\tilde{\delta}}(O(t)))^c$. Moreover, we assume that $| dist(\gamma_i, \gamma_j)|>2\delta_i$ for every $i\ne j$. 
%Here $i+1, i-1$ are considered modulus 3 and $c$ denotes the complement.

Now define $$\delta=\min_{i=1,2,3} \delta_i.$$

\item Consider the following sets: $$D_{ii}(t)=\{x\in \Omega: d_i(x,t)\leq \delta\}$$
$$D_{ii+1}(t)=\left\{x\in \Omega: d_i(x,t)\geq \frac{\delta}{2}, \quad d_{i+1}(x,t)\leq -\frac{\delta}{2} \hbox{ and }d_{ii+1}\geq\frac{\delta}{2} \right\}\subset S_{ii+1}.$$
Define a partition of unity $\xi^{ext}_{ij}(x,t)$ associated to these sets, that is define functions 
$0\leq \xi^{ext}_{ij}\leq 1$ such that $j\in \{i,i+1\}$, $supp \ \xi^{ext}_{ij}\subset D_{ij}$ (where $supp$ denotes the support) and for every $x\in\bigcup_{i,j}D_{ij}$ holds $\sum_{i,j}\xi^{ext}_{ij}(x,t)=1$.
It is easy to see that is possible to define these functions so that
$$|\nabla \xi^{ext}_{ij}(x,t)|\leq \frac{C}{\delta}\max_{i,j}|\nabla d_{i}(x)|$$
for $x$ such that $\frac{\delta}{2}\leq |d_i|(x,t)\leq \delta $ or $\frac{\delta}{2}\leq | d_j|(x,t)\leq \delta$ and
$$|\nabla \xi^{ext}_{ij}(x,t)|=0$$  when   $|d_i|(x,t)\geq \delta $ and  $|d_j|(x,t)\geq \delta$ or
 when $i=j$ and $|d_i|(x,t)\leq \frac{\delta}{2} $.

 \item  We define the following boundary condition:

\begin{align}
\phi_\epsilon(x,t)=
\sum_{i=1}^3 &\left( \xi^{ext}_{ii}\left(x,t\right)
\zeta_{ii+1}\left(
\frac{d_i(x,t)}{\epsilon}\right) \xi^{ext}_{ii+1}\left(x,t\right)
c_i  \right). \label{defbc}\end{align} 

Notice that this function it is well defined only in $\bigcup_{i,j}D_{ij}$. In the next step we choose an appropriate cut-off function   to extend it to the whole domain $\Omega$.

\item Let  $r(x,t)=|x-O(t)|$ and define

$$\phi^\eta_\epsilon(x,t)\equiv \left(1-\eta_2\left(\frac{r(x,t)}{2\epsilon}+1-\frac{\tilde{\delta}}{2\epsilon}\right)\right)\phi_\epsilon(x,t).$$

\item Similarly inside the $B_{\tilde{\delta}}(O(t))$ we define

$$\tilde{\phi}_\epsilon(x,t)=
\sum_{i=1}^3\left( \eta_{2i}(\theta-\theta(t))
\zeta_{ii+1}\left(
\frac{d_i(x,t)}{\epsilon}\right)
\frac{}{}+\eta_{2i-1}(\theta-\theta(t))c_i\right).$$  

\item Let $\frac{1}{2}<\rho<1$. We extend $\tilde{\phi}$ to the whole domain by
$$\tilde{\phi}^\eta_\epsilon(x,t) \equiv \left(1-\eta_1\left(\frac{x-O(t)}{\epsilon^\rho}\right)\right) \tilde{\phi}_\epsilon(x,t)$$

\item Now we let
\begin{align}\tilde{v}_\epsilon(x,t)=&\tilde{\phi}^\eta_\epsilon(x,t)+\eta_1\left(\frac{x-O(t)}{\epsilon^\rho}\right)u_*\left( \frac{R_{\theta(t)}(x-O(t))}{\epsilon}\right).\end{align}

\item Finally we define
\begin{align} v_\epsilon(x,t)=&\phi^\eta_\epsilon(x,t)
+\eta_2\left(\frac{r(x,t)}{2\epsilon}+1-\frac{\tilde{\delta}}{2\epsilon}\right)
\tilde{v}_\epsilon(x,t).\label{defv}\end{align}

 We let
\be \psi_\epsilon(x)=v_\epsilon(x,0) \label{defci}\ee
to be the initial condition in \equ{ci}.

%\item To simplify the notation in what follows, we denote

\end{itemize}

\medskip

Then we prove:

\begin{teo}\label{mainteo}
Suppose that we have a triod $\calT=\{\gamma^i\}$ satisfying \equ{triod} in a compact domain $\Omega$
 Additionally, we assume that at  the meeting point $O(t)$ the angles formed by the tangents are prescribed and fixed, that is
there are
$0<\alpha_1<\alpha_2<\alpha_3<2\pi$ such that
$\theta^i(t)=\theta^{i-1}(t)+\alpha_i$ for every $t>0$.
 Assume also that $\calT$ is well defined for every $t\in [0,T]$, that the meeting  point 
 satisfies $\min_{t\in[0,T]}dist(O(t),\partial \Omega)>0$, $ |O'(t)|\leq C$ and $|\theta'(t)|\leq C$.
If there is a potential  $W$ satisfying conditions (W1)-(W4) such that \equ{condang} is satisfied (for the $\alpha_i$ defined above) then the unique   solution $u_\epsilon$ to \equ{laeq}-\equ{ci}-\equ{bc} (where $\phi_\epsilon$ and $\psi_\epsilon$ are given by \equ{defci} and \equ{defbc} respectively)
that satisfies
\be \lim_{\epsilon \to 0}\sup_{\Omega\times[0,T]}|u_\epsilon-v_\epsilon|(x,t)=0, \ee
for $v_\epsilon$ is given by \equ{defv}.

\end{teo}

\begin{rem}
We will assume that at the boundary point the mean curvature flow equation is satisfied, i.e.
\be \frac{\partial}{\partial t}P^i(t)=\frac{\gamma^i_{xx}}{|\gamma^i_x|^2}(1,t).\label{cbg}\ee

Notice that the computations in \cite{moybycm} show that this  holds when we consider fixed end-points. 
%For non- compacts domains $\Omega$ condition \equ{bc} should be understood 
%as there is a given $\phi_\epsilon$ such that
%$$\lim_{x\in \Omega, |x|\to \infty}|u_\epsilon(x,t)-\phi_\epsilon(x,t)|=0.$$ 
\end{rem}

We also prove the 2 following corollaries:

\begin{cor}\label{cor1}
Suppose that there are two triods $\calT_1$ and $\calT_2$ satisfying the conditions of Theorem \ref{mainteo}. Moreover, assume that for positive times 
they satisfy the same meeting condition at $O^i(t)$ (i.e. the prescribed angles $\theta_i$ are the same). Assume that  
$$\calT_1(\cdot,0)=\calT_2(\cdot,0)$$
as sets.
% (i.e one is a reparametrization of the other).
Then for every $t>0$ holds
$$\calT_1(\cdot,t)=\calT_2(\cdot,t)$$
as sets. That is, there is a unique geometric solution to \equ{triod}. 
\end{cor}

\begin{cor}\label{cor2}
%Notice that the angle condition in \equ{triod} it is only required for $t>0$. Hence we have:

%\begin{teo}
Consider a triod $\calT=\{\gamma^i\}_{i=1}^3$ satisfying \equ{triod} and the conditions in Theorem \ref{mainteo}. Suppose in addition that that for every $t>0$ $\calT$ satisfies 
$$\sum_{i=1}^3 \frac{\gamma_\lambda^i(0,t)}{|\gamma_\lambda^i(0,t)|}=0.$$
That is, the curves meet at $120^0$ for every $t>0$ (but not necessarily at $t=0$).
Suppose the $O(0)=\gamma^i(0,0)=0$ and that there is a constant $C$, such that for each $i$
$$\sup_i|k_i(\lambda,t)|\leq \frac{C}{\sqrt{t}},$$ 
where $k_i(\cdot,t)$ is the curvature of the curve $\gamma^i(\cdot,t)$.
Let $\beta_n \to 0 $ be a sequence of positive real numbers.
Then the sequence of triods defined by
$\calT_n=\left\{\frac{1}{\beta_n}\gamma^i(\cdot,\beta_n^2 t)\right\}$ converges uniformly in compact sets to the self-similar solution described in \cite{selsimos} with initial condition $\left\{ \lambda \frac{\gamma_\lambda^i(0,t)}{|\gamma_\lambda^i(0,t)|}, \quad \lambda \in \rr_+ \right\}_{i=1}^3$. 
%\end{teo}

\end{cor}
In order to prove Theorem \ref{mainteo} we need 
the following functional:
%Let $\calQ=(0,1]\times(0,1]\times [0,1]$.
%Define for  ${\vec{q}}=
%(\epsilon,\sigma, \alpha)\in \calQ$ such that $\sigma\leq \epsilon^{1-\alpha}$ the function
%  $$v_{\vec{q}}(y)=\eta_{\alpha}(y) \phi (y)+(1-\eta_{\alpha}(y)) v_\sigma (y).$$ 
%For $\vec{q}$ as above consider the functional 
 $$F_ {\epsilon}(h,\psi_\epsilon)=-\int_0^t \int_{\Omega}\calH_{\Omega}(x,y,t-s)\left(\frac{\nabla_u W(h+\phi^\eta_{\epsilon})}{\epsilon^2}+P \phi^\eta_\epsilon\right)(y,s)dyds$$
$$+ \int_{\Omega}\calH_{\Omega}(x,y,t)(\psi_\epsilon (y)-\phi^\eta_\epsilon(y,0))dy,$$
where $\calH_{\Omega}$ denotes the heat kernel in  $\Omega$.
 
\begin{rem}\label{uh} Notice that fixed points of this functional are solutions to the equation
\begin{align}
\frac{\partial h_\epsilon}{\partial t}-\Delta h_\epsilon
+\frac{\nabla_u W(h_\epsilon+\phi^\eta_\epsilon)}{2\epsilon^2}&=-P \phi^\eta_\epsilon\quad \hbox{ in }\Omega\label{ecpar}\\
h_\epsilon(x,t)&=0 \hbox{ on }\partial \Omega\label{cbpar}\\
h_\epsilon(x,0)&=\psi(x)-\phi^\eta_\epsilon(x,0).\label{cipar}
\end{align}

In particular, defining $u_\epsilon(x,t)=h_\epsilon(x,t)+\phi^\eta(x,t)$ we have $u_\epsilon(x,t)$ satisfies \equ{laeq}-\equ{ci}-\equ{bc}.
\end{rem}

The main tool that we will use to prove  Theorem \ref{mainteo} is Lemma 4.1 in \cite{stationary}. We restate it here without proof.

\begin{lem}\label{cotaprinc}

Fix $K>0$. Consider the sequences of continuous functions $\psi_n,
w_n $ satisfying $\sup |\psi_n|,\ \sup|w_n|\leq K$. Let
 $\epsilon_n\to 0$ and $T_n>0$. 
Assume in addition that for every $0<\epsilon<1$ holds
$\sup_{x\in \Omega, t\in [0,T]}|h_\epsilon|(x,t)\leq K$
Then for each $\psi_n,\epsilon_n$ the functional $F_{\epsilon_n}$ has a unique fixed point $h_{\epsilon_n}$ and holds either
\begin{enumerate}
\item \label{mejorcaso}$\lim_{n \to \infty}\sup_{\Omega\times [0,T_n]}|w_n-h_{\epsilon_n}|\to 0$,
or
\item \label{peorcaso} there is a constant $C$, independent of ${\epsilon}_n$ and $T_n$
such that 
$$\sup_{\Omega\times [0,T_n]}|w_n-h_{\epsilon_n}|\leq C \sup_{\Omega\times [0,T_n]}
|F_{\epsilon_n}(w_n,\psi_n)-w_n|.$$

\end{enumerate}
\end{lem}

%\begin{rem}\label{uh}
%Notice that defining $u_\epsilon(x,t)=h_\epsilon(x,t)+\phi(x,t)$ we have that $u_\epsilon(x,t)$ satisfies \equ{laeq}-\equ{ci}-\equ{cb}. 
%\end{rem}

From the proof of Lemma 4.1 in \cite{stationary} we have the following corollary:

\begin{cor}\label{corcotprinc}
Consider $w_n$ and $h_{\epsilon_n}$ as in Lemma \ref{cotaprinc}.
Suppose that there is $l>0$ such that
$$\frac{|F_{\epsilon_n}(w_n,\psi_n)-w_n|}{\epsilon_n^l}\to 0 \hbox{ as } \epsilon\to 0,$$
then
$$\frac{|h_{\epsilon_n}-w_n|}{\epsilon_n^l}\to 0 \hbox{ as } \epsilon\to 0.$$

\end{cor}

Regarding a priori bounds, existence and uniqueness of solutions, we note that
 Theorems 4.2 and 6.2 in \cite{stationary} can be easily extended to our setting for any compact domain $\Omega$. That is, there exists a unique solution to \equ{ecpar}, \equ{cipar} and \equ{cbpar} that satisfies $|h_{\epsilon}(x,t)|\leq C$, where $C$ depends only on $W$, $\sup |\phi_\epsilon|$ and $\sup| \psi_\epsilon|$ (in particular can be chosen independent of $\epsilon$ if $\sup |\phi_\epsilon|$ and $\sup| \psi_\epsilon|$  are bounded independently of $\epsilon$).

We would like to point out that the computations in the proof of Theorem \ref{mainteo}  are similar to the ones in \cite{stationary} and \cite{papertesis}. We refer the reader to these papers for further details in the calculations.

%\section{ Existence and uniqueness in non-compact Domains}

%In order to prove existence for non-compact domains we consider the following sequence of solutions:
%\begin{align}\frac{\partial u^R_\epsilon}{\partial t}-\Delta u^R_\epsilon+ \frac{\nabla_uW(u^R_\epsilon)}{\epsilon^2}&=0 \hbox{ for } x\in \Omega \bigcap B_R\label{laeqr}\\
% u^R_\epsilon(x,0)&=\psi^R_\epsilon(x),\label{cir} \\
%u_\epsilon|_{\partial \Omega\bigcap B_R}&=\phi_\epsilon (x,t), \label{cbr}\end{align}
%where
%$\psi_\epsilon^R(x)=\chi\left(\right)\psi_\epsilon+\left(1-\chi\left(\right)\right)\phi_\epsilon(x,0)$
\section{Proof of Theorem \ref{mainteo}}\label{pfth}

\begin{proof}

To prove Theorem \ref{mainteo} we use Lemma \ref{cotaprinc} with the sequence
$$w_\epsilon(x,t)=\eta_2\left(\frac{r(x,t)}{2\epsilon}+1-\frac{\tilde{\delta}}{2\epsilon}\right) 
\tilde{v}_\epsilon(x,t).$$
Notice that  $w_\epsilon(x,t)+\phi^\eta_\epsilon(x,t)=v_\epsilon(x,t)$. Moreover, without loss of generality we can assume that in  \equ{gbdeq} we chose $\tilde{\delta}<\min_{[0,T]}dist(O(t),\partial \Omega)$, therefore
 $w_\epsilon|_{\partial \Omega}=0$ and
$$w_\epsilon(x,t)=\int_0^t \int_{\Omega}\calH_{\Omega}(x,y,t-s)
P w_\epsilon(y,s) dyds
+ \int_{\Omega}\calH_{\Omega}(x,y,t)w_\epsilon (y,0)dy.$$ Which implies  (recall that $\psi_\epsilon(y)=v_\epsilon(y,0)$),
\begin{align}(F_\epsilon(w_\epsilon,\psi_\epsilon)-w_\epsilon)(x,t)=&\int_0^t \int_{\Omega}\calH_{\Omega}(x,y,t-s)\left(-\frac{\nabla_u W(w_\epsilon+\phi^\eta_{\epsilon})}{\epsilon^2}(y,s)-P (\phi^\eta_\epsilon+w_\epsilon)(y,s)\right)dyds\notag \\
&+ \int_{\Omega}\calH_{\Omega}(x,y,t)(\psi_\epsilon (y)-\phi^\eta_\epsilon(y)-w_\epsilon(y,0))dy\notag \\
=&\int_0^t \int_{\Omega}\calH_{\Omega}(x,y,t-s)\left(-\frac{\nabla_u W(v_\epsilon)}{\epsilon^2}-P v_\epsilon\right)(y,s)dyds.
\label{eqvint}\end{align}
Consider a sequence of $\epsilon_n\to 0$. Using Remark \ref{uh} and
 Lemma \ref{cotaprinc} we have that either
\begin{enumerate}
\item $\lim_{n \to \infty}\sup_{\Omega\times [0,T]}|v_{\epsilon_n}-u_{\epsilon_n}|=
\lim_{n \to \infty}\sup_{\Omega\times [0,T]}|w_{\epsilon_n}-h_{\epsilon_n}|\to 0$,
or
\item \label{peorcaso} there is a constant $C$, independent of ${\epsilon}_n$ and $T$
such that 
$$\lim_{n \to \infty}\sup_{\Omega\times [0,T]}|v_{\epsilon_n}-u_{\epsilon_n}|=
\lim_{n \to \infty}\sup_{\Omega\times [0,T]}|w_{\epsilon_n}-h_{\epsilon_n}|
\leq C \sup_{\Omega\times [0,T]}
|F_{\epsilon_n}(w_n,\psi_{\epsilon_n})-w_n|.$$
\end{enumerate}

Suppose that we are in the second case. We will show that  $\sup_{\Omega\times [0,T]}
|F_{\epsilon_n}(w_n,\psi_{\epsilon_n})-w_n|\to 0$ concluding the result.

Recalling equation \equ{eqvint},  we  
compute $Pv_\epsilon+\frac{\nabla _u W(v_\epsilon)}{\epsilon^2}$.  By the definition of the function $v_\epsilon$,  we notice that for each fixed time $t$  the cut-off functions $\eta_i$ divide the space into three interior regions   (namely $B_{\frac{\epsilon^\rho}{2}}(O(t))$,  
$B_{\tilde{\delta}- \epsilon}(O(t))\setminus B_{\epsilon^\rho}(O(t))$ and $\Omega\setminus B_{\tilde{\delta}}$)
and into  2 transition regions ($B_{\epsilon^\rho}(O(t))\setminus  B_{\frac{\epsilon^\rho}{2}}(O(t))$ and $B_{\tilde{\delta}}(O(t))\setminus B_{\tilde{\delta}- \epsilon}(O(t))$). We compute separately in each of them.

\medskip

\noindent{ \large{\bf Interior Regions:}}
\smallskip
\begin{itemize}
\item{\em In $B_{\frac{\epsilon^\rho}{2}}(O(t))$}:

Inside $B_{\frac{\epsilon^\rho}{2}}(O(t))$ we have $v_\epsilon(x,t)=u_*\left(\frac{R_{\theta(t)}(x-O(t))}{\epsilon}\right) $. For simplicity, in the computation that follows we are going to omit the argument of the function $u_*$ and its derivatives (which will always be $\frac{R_{\theta(t)}(x-O(t))}{\epsilon}$). Hence,

\begin{align*}\left(Pv_\epsilon+\frac{\nabla_u W\left( v_\epsilon\right)}{\epsilon^2}\right)(x,t)=& Ju_*\frac{\theta' R_{\theta(t)+\frac{\pi}{2}}(x-O(t))-R_{\theta(t)}O'(t)}{\epsilon} - \frac{\Delta u_*}{\epsilon^2}
+ \frac{\nabla_u W(u_*)}{\epsilon^2}\\
 =
& Ju_*\frac{\theta' R_{\theta(t)+\frac{\pi}{2}}(x-O(t))-R_{\theta(t)}O'(t)}{\epsilon}  \end{align*}
Here $Ju_*$ denotes the Jacobian matrix of $u_*$.
By our hypotheses on $\theta'$ and $O'(t)$ we have that there is a constant $C$ such that
$$\left|Pv_\epsilon+\frac{\nabla_u W\left( v_\epsilon\right)}{\epsilon^2}\right|(x,t)\leq \frac{C}{\epsilon}.$$
 This implies
\begin{align*}\left|\int_0^t \int_{B_{\frac{\epsilon^\rho}{2}}(O(s))}\calH_{\Omega}(x,y,t-s)\left(\frac{}{}\right.\right.& \left.\left.\frac{\nabla_u W(v_\epsilon)}{\epsilon^2}(y,s)
+P v_\epsilon(y,s)\right)dyds\right|\\
\leq &\frac{C}{\epsilon}\int_0^t \int_{B_{\frac{\epsilon^\rho}{2}}(O(s))} \calH_{\Omega}(x,y,t-s)dy ds\end{align*}

Recall that we assume that $dist(O(t),\partial \Omega)>\tilde{\delta}$ uniformly on $t$. Therefore there is another constant $C$ such that for every $x\in \Omega$, $y\in B_{\frac{\epsilon^\rho}{2}}(O(s))$, $s\in [0,T]$ holds
$\left|\calH_{\Omega}(x,y,t-s)\right|\leq \frac{Ce^{-\frac{|x-y|^2}{t-s}}}{t-s}$. Thus,
\begin{align}\left|\int_0^t \int_{B_{\frac{\epsilon^\rho}{2}}(O(s))}\calH_{\Omega}(x,y,t-s)\left(\frac{}{}\right.\right.& \left.\left.\frac{\nabla_u W(v_\epsilon)}{\epsilon^2}(y,s)
+P v_\epsilon(y,s)\right)dyds\right|\notag \\
\leq & \frac{C}{\epsilon}\int_0^t \int_{B_{\frac{\epsilon^\rho}{2}}(O(s))} \frac{e^{-\frac{|x-y|^2}{t-s}}}{t-s} dy ds\notag \\
\leq &\frac{C}{\epsilon}\left(\int_0^{t-\epsilon^m} \int_{B_{\frac{\epsilon^\rho}{2}}(O(s))} \frac{1}{t-s} dy ds+\int_{t-\epsilon^m}^t \int_{\rr^2} \frac{e^{-\frac{|x-y|^2}{t-s}}}{t-s} dy ds \right)\notag \\
\leq &\frac{C}{\epsilon}\left(\int_0^{t-\epsilon^m} \frac{\pi \epsilon^{2\rho}}{4(t-s)}ds+ \int_{t-\epsilon^m}^t ds \right)\notag \\ \leq &\frac{C}{\epsilon}\left((\ln t-m\ln\epsilon)\frac{\pi \epsilon^{2\rho}}{4}
+ \epsilon^m\right)\notag \\ \leq &C\left((\ln T-m\ln\epsilon)\frac{\pi \epsilon^{2\rho-1}}{4}
+ \epsilon^{m-1}\right)\to 0 \hbox{ as } \epsilon \to 0. \label{bd1}\end{align}
The convergence as $\epsilon \to 0$ holds since $\rho>\frac{1}{2}.$
This completes for the moment the bound needed in $B_{\frac{\epsilon^\rho}{2}}(O(t))$.

\medskip 

Since the computations in the 2 remaining interior regions ($B_{\epsilon^\rho}(O(t))\setminus
B_{\tilde{\delta}- \epsilon}(O(t))$ and $\Omega\setminus B_{\tilde{\delta}}(t)$) are similar to each other,   we will only present the one in $B_{\tilde{\delta}- \epsilon}(O(t))\setminus B_{\epsilon^\rho}(O(t))$ in detail and  point out the necessary modifications in $\Omega\setminus B_{\tilde{\delta}}(t)$.

\smallskip

\item {\em In $B_{\tilde{\delta}- \epsilon}(O(t))\setminus B_{\epsilon^\rho}(O(t))$:}

It  holds that $v_\epsilon(x,t)=\tilde{\phi}(x,t)$ and $\eta_1\equiv 1$, then
\begin{align}\left|Pv_\epsilon+\frac{\nabla_u W\left( v_\epsilon\right)}{\epsilon^2}\right|(x,t)= &\left| \sum_{i=1}^3\left(P\left( \xi^{ext}_{2i}(\theta-\theta(t))\right)
\zeta_{ii+1}\left(
\frac{d_i(x,t)}{\epsilon}-c_i\right)\right.\right.\notag \\
&+\xi^{ext}_{2i}(\theta-\theta(t))
P\left(\zeta_{ii+1}\left(
\frac{d_i(x,t)}{\epsilon}\right)\right)\notag 
\\
&\left. -2\nabla\xi^{ext}_{2i}(\theta-\theta(t)) \cdot J\zeta_{ii+1}\left(
\frac{d_i(x,t)}{\epsilon}\right) \right)
\notag \\ &  \left.+\frac{\nabla_u W\left( v_\epsilon\right)}{\epsilon^2}\right|.
\label{eqout}\end{align} 

where $J\zeta_{ii+1}=\left(\frac{\partial}{\partial x_m}\zeta^l_{ii+1}\right)_{lm}$ is the Jacobian of $\zeta_{ii+1}$ 
and 
the $j$-th component of the vector 
$\nabla\xi^{ext}_{2i}(\theta-\theta(t)) \cdot J\zeta_{ii+1}$ is defined by
$\nabla\xi^{ext}_{2i}(\theta-\theta(t)) \cdot (J\zeta_{ii+1})_{j\cdot}$, where $\cdot$ is the standard dot product and $ (J\zeta_{ii+1})_{j\cdot}$ denotes the $j$-th row.

 For each fixed time $t$ we can 
write every $x \in B_{\tilde{\delta}- \epsilon}(O(t))\setminus 
B_{\epsilon^\rho}(O(t))$ as $$x=O(t)+r(cos \theta, sin \theta).$$ 
By the definition of the functions $\xi^{ext}_{i}$,
within   $B_{\tilde{\delta}- \epsilon}(O(t))\setminus B_{\epsilon^\rho}(O(t))$ we can distinguish 3 spatial types of regions that depend on $\theta$: near an interfaces (i.e.  $x$ is near an interface if
there is an $i_0$ such that
$\xi^{ext}_{2i_0}(\theta-\theta(t))=1$ and $\xi^{ext}_{j}(\theta-\theta(t))=0$ for $j\ne 2i_0 $), away from the interfaces (i.e.  $x$ is away of the interfaces if there is an $i_0$ such that
$\xi^{ext}_{2i_0-1}(\theta-\theta(t))=1$ and  $\xi^{ext}_{j}(\theta-\theta(t))=0$ for $j\ne 2i_0-1 $) and the transition regions (which correspond to  $x$ such that there is
an  $i_0$ such that $\xi^{ext}_{j}(\theta-\theta(t))=0$ for $j\ne 2i_0-1, 2i_0$,  
$\xi^{ext}_{2i_0-1}(\theta-\theta(t)),\xi^{ext}_{2i_0}(\theta-\theta(t))\ne 1$
 and
$\xi^{ext}_{2i_0-1}(\theta-\theta(t))+\xi^{ext}_{2i_0}(\theta-\theta(t))=1$).

In each of these regions several of the terms in equation \equ{eqout} cancel. Therefore we compute separately in each of them. Let us consider the three possible cases:

\medskip

\begin{enumerate}
\item {\bf Near an interface:}

 If there is an $i_0$ such that
$\xi^{ext}_{2i_0}(\theta-\theta(t))=1$ and $\xi^{ext}_{j}(\theta-\theta(t))=0$ for $j\ne 2i_0 $,
then  $\tilde{\phi}(x,t)=\zeta_{i_0 i_0+1}\left(\frac{d_{i_0}(x,t)}{\epsilon}\right)$ and equation \equ{eqout} reduces to
\begin{align}\left|P\zeta_{i_0i_0+1}+\frac{\nabla_u W\left(\zeta_{i_0i_0+1}\right)}{\epsilon^2}\right|=&\left|
\left(
\frac{(d_{i_0}(x,t))_t-\Delta d_{i_0}(x,t)}{\epsilon}\right) \zeta'_{i_0i_0+1}\right.\notag \\
&\left.-\frac{\zeta''_{i_0i_0+1}}{\epsilon^2}
|\nabla d_{i_0}(x,t)|^2+ \frac{\nabla_u W\left(\zeta_{i_0i_0+1}\right)}{\epsilon^2}\right|.\label{eqoutin}
\end{align}
In this last equation we omitted the argument of the function $\zeta_{i_0i_0+1}$ and its derivatives, but these arguments should always be $\frac{d_{i_0}(x,t)}{\epsilon}$. Since the interfaces $\gamma_{i_0}$ evolve under curvature flow we have that the distance function satisfy:
\be (d_{i_0})_t-\Delta d_{i_0}=\frac{k_{i_0}^2(\lambda,t)d_{i_0}}{1+k_{i_0}(\lambda,t)d_{i_0}} \hbox{ and }|\nabla d_{i_0}|=1 \label{ecdist}\ee

where $k_{i_0}(\lambda,t)$ is the curvature of $\gamma_{i_0}(\lambda,t)$ at  the point where the distance $d_{i_0}(x,t)$
is attained at time $t$ (for details on this computations see \cite{calofla} or \cite{tesis} for example).  Combining \equ{eqoutin}, \equ{ecdist} and the exponential decay
 of $\zeta_{i_0i_0+1}$ and its derivatives we have that
\equ{eqout} in this region equals to
\begin{align}\left|Pv_\epsilon+\frac{\nabla_u W\left( v_\epsilon\right)}{\epsilon^2}\right|=&\left| \frac{k_{i_0}^2(\lambda,t)d_{i_0}}{\epsilon(1+k_{i_0}(\lambda,t)d_{i_0})}\zeta'_{i_0i_0+1}\right|\notag\\
\leq &
C\left| \frac{k_{i_0}^2(\lambda,t)d_{i_0}}{\epsilon(1+k_{i_0}(\lambda,t)d_{i_0})}\right| e^{-c\frac{d_{i_0}(x,t)}{\epsilon}}\label{ec1int}\end{align}

\smallskip

\item { \bf Away from the interfaces:} 

If
there is an $i_0$ such that
$\xi^{ext}_{2i_0-1}(\theta-\theta(t))=1$ and  $\xi^{ext}_{j}(\theta-\theta(t))=0$ for $j\ne 2i_0-1 $,
then  $\tilde{\phi}(x,t)=c_{i_0}$
and equation \equ{eqout} is identically 0:
\be \left|Pv_\epsilon+\frac{\nabla_u W\left( v_\epsilon\right)}{\epsilon^2}\right|(x,t)=0 \label{ecint2}\ee

\smallskip

\item {\bf  Transition Regions:}

If there is
an  $i_0$ such that $\xi^{ext}_{j}(\theta-\theta(t))=0$ for $j\ne 2i_0-1, 2i_0$, $\xi^{ext}_{2i_0-1}(\theta-\theta(t)),\xi^{ext}_{2i_0}(\theta-\theta(t))\ne1$ and
$\xi^{ext}_{2i_0-1}(\theta-\theta(t))+\xi^{ext}_{2i_0}(\theta-\theta(t))=1$,

then, by the definition of $\xi^{ext}_j$ and $\theta_j$, we have  that $|\theta-\theta_j|>\delta_{int}$ for every  $j=1,2,3$. 

Let $$y=O(t)+r(y)(\cos \theta(y), \sin \theta(y)) \in \gamma_j(\cdot,t).$$ Then it holds  $|\theta (y)-\theta_j|\leq \frac{\delta_{int}}{2}$.
Therefore, $|\theta (y)-\theta|\geq \frac{\delta_{int}}{2}$ for every $y=\gamma_j(\lambda,t)$. 
This implies $dist(x,y)\geq |x-O(t)|\cos \frac{\delta_{int}}{2}$. Since $x\in B_{\tilde{\delta}- \epsilon}(O(t))\setminus 
B_{\epsilon^\rho}(O(t))$, we conclude  
\be d_j(x,t)>\epsilon^\rho \cos \frac{\delta_{int}}{2}.\label{bdd}\ee
By writing 
$\frac{\nabla_u W\left( v_\epsilon \right)}{\epsilon^2}=
\frac{\nabla_u W\left( v_\epsilon \right)-\nabla_u W(c_i)}{\epsilon^2} $
we have
$$\left|\frac{\nabla_u W\left( v_\epsilon \right)}{\epsilon^2}\right| \leq
C \frac{ \left| v_\epsilon -c_i\right|}{\epsilon^2}, $$
where $C$ depends on the second derivatives of $W$ and the uniform bounds of $v_\epsilon$.
Now using  \equ{ecdist}, \equ{bdd} and the exponential decay of $\zeta_{i_0i_{0+1}}$ we have the following bound for \equ{eqout}:
\be  \left|Pv_\epsilon+\frac{\nabla_u W\left( v_\epsilon \right)}{\epsilon^2}\right|(x,t)\leq C\frac{e^{-c\frac{d_{i_0}(x,t)}{\epsilon}}}{\epsilon^2}\leq 
 C\frac{e^{-c\epsilon^{\rho-1}}}{\epsilon^2},\label{ecint3}\ee
for some constants $C,c$ that depend on $\delta_{int}$ and $\zeta_{ij}$.

\end{enumerate}

\medskip 

From equations \equ{ec1int}, \equ{ecint2} and \equ{ecint3} we can see that
$$\left|Pv_\epsilon+\frac{\nabla_u W\left( v_\epsilon\right)}{\epsilon^2}\right|\leq C \max \left\{\sup_i | k_i|,
 C\frac{e^{-c\epsilon^{\rho-1}}}{\epsilon^2}\right\} $$
and
$$\left|Pv_\epsilon+\frac{\nabla_u W\left( v_\epsilon\right)}{\epsilon^2}\right|(x,t) \to 0 \hbox{ as } \epsilon \to 0$$  for every $ x \in  B_{\tilde{\delta}- \epsilon}(O(t))\setminus B_{\epsilon^\rho}(O(t)), t>0.$

Therefore, if $\sup_i|k_i|\leq C$ for every $t>0$, we have
\begin{align}\int_0^t &\int_{ B_{\tilde{\delta}- \epsilon}(O(t))\setminus B_{\epsilon^\rho}(O(t))}
 \calH_{B_1}(x,y,t-s) \left|Pv_\epsilon+\frac{\nabla_u W\left( v_\epsilon\right)}{\epsilon^2}\right|(y,s)
dy ds \notag \\ = & \int_{t-\delta}^t \int_{ B_{\tilde{\delta}- \epsilon}(O(t))\setminus B_{\epsilon^\rho}(O(t))}
\calH_{B_1}(x,y,t-s) \left|Pv_\epsilon+\frac{\nabla_u W\left( v_\epsilon\right)}{\epsilon^2}\right|(y,s)
dy ds  \notag \\ &+\int_0^{t-\delta} \int_{ B_{\tilde{\delta}- \epsilon}(O(t))\setminus B_{\epsilon^\rho}(O(t))}
\calH_{B_1}(x,y,t-s) \left|Pv_\epsilon+\frac{\nabla_u W\left( v_\epsilon\right)}{\epsilon^2}\right|(y,s)
dy ds \notag \\
\leq & C\delta+ \frac{1}{\delta}\int_0^t \int_{ B_{\tilde{\delta}- \epsilon}(O(t))\setminus B_{\epsilon^\rho}(O(t))}
\left|Pv_\epsilon+\frac{\nabla_u W\left( v_\epsilon\right)}{\epsilon^2}\right|(y,s)
dy ds  \to C\delta \hbox{ as }\epsilon \to 0. \label{bd2}\end{align}

\begin{rem}\label{difbal}
It is easy to see that for any $r(\epsilon)\leq \tilde{\delta}$ and  $|x-O(t)|\geq r(\epsilon)$ equation 
\equ{bdd} can be replaced by
\be d_j(x,t)>r(\epsilon) \cos \frac{\delta_{int}}{2}.\label{bdd2}\ee
for $x\in B_{\tilde{\delta}- \epsilon}(O(t))\setminus 
B_{r(\epsilon)}(O(t))$. In particular, for $r(\epsilon)=K\epsilon^\rho$
the estimates above hold for a different constant $C$, that depends on $K$.

\end{rem}
\smallskip

\begin{rem} \label{unifbounds}
Using  \equ{ec1int},  \equ{ecint2} and  \equ{ecint3} it is easy to see that the bound above can be computed more precisely. Namely, by separating the domain into $|d_i|\leq \epsilon^\rho$ and $|d_i|\geq \epsilon^\rho$ for any $\rho<1$ it is easy to compute that we have
$$\int_0^t \int_{ B_{\tilde{\delta}- \epsilon}(O(t))\setminus B_{\epsilon^\rho}(O(t))}
\left|Pv_\epsilon+\frac{\nabla_u W\left( v_\epsilon\right)}{\epsilon^2}\right|(y,s)
dy ds \leq Ce^{-c\epsilon^{\rho-1}}+C\epsilon^{2\rho},$$
where the constants depend only on uniform bounds of the curvatures $k_i$.

\end{rem}

\begin{rem}\label{remss}
 The computation above carries over in a similar way when  $\sup_i|k_i|\leq f(t)$ where $f(t)$ is an integrable function of $t$. In this case the bounds above will depend on $\int_0^t|k_i|$. 

A particular and important example are the self-similar solutions computed in \cite{selsimos}, where
$|k^i|\leq \frac{C}{\sqrt{t}}$. The computations above gives:
\begin{align*}\int_0^t \int_{ B_{\tilde{\delta}- \epsilon}(O(t))\setminus B_{\epsilon^\rho}(O(t))}
 & \calH_{B_1}(x,y,t-s)  \left|Pv_\epsilon+\frac{\nabla_u W\left( v_\epsilon\right)}{\epsilon^2}\right|(y,s)
dy ds \\ \leq & C\int_{t-\delta}^\delta\frac{1}{\sqrt{s}}ds+ \frac{1}{\delta}\int_0^t \int_{ B_{\tilde{\delta}- \epsilon}(O(t))\setminus B_{\epsilon^\rho}(O(t))}
\left|Pv_\epsilon+\frac{\nabla_u W\left( v_\epsilon\right)}{\epsilon^2}\right|(y,s)
\\  \leq &
C(\sqrt{t}-\sqrt{t-\delta})+C\frac{\epsilon^{2\rho}  }{\delta}\int_{0}^t\frac{1}{\sqrt{s}}ds \\
\leq &C\delta+\frac{C}{\delta}\sqrt{T}\epsilon^{2\rho},\end{align*}
where the constants are independent of $\epsilon$ and depend linearly  on uniform bounds of $\sqrt{t}k_i$.

\end{rem}

\item{\em In $\Omega\setminus B_{\tilde{\delta}}$:
}

In this set $v(x,t)=\phi(x,t)$. As in the previous case, the function $\chi_{ij}$  divide
$\Omega\setminus B_{\tilde{\delta}}$ into regions like the ones described above: close to the interface, away form the interface and transition regions. The bounds in the different sets are analogous to the ones in $B_{\tilde{\delta}- \epsilon}(O(t))\setminus B_{\epsilon^\rho}(O(t)) $. We find that
\begin{align}\left|Pv_\epsilon+\frac{\nabla_u W\left( v_\epsilon\right)}{\epsilon^2}\right|=&\left| \sum_{i=1}^3\left(P\left( \xi^{ext}_{2i}(\theta-\theta(t))\right)
\zeta_{ii+1}\left(
\frac{d_i(x,t)}{\epsilon}-c_i\right)\right.\right.\notag \\
&+\xi^{ext}_{2i}(\theta-\theta(t))
P\left(\zeta_{ii+1}\left(
\frac{d_i(x,t)}{\epsilon}\right)\right)\notag 
\\
&\left. -2\nabla\xi^{ext}_{2i}(\theta-\theta(t)) \cdot J\zeta_{ii+1}\left(
\frac{d_i(x,t)}{\epsilon}\right) \right)
\notag \\ &  \left.+\frac{\nabla_u W\left( v_\epsilon\right)}{\epsilon^2}\right|
\notag\\
\leq &
C\max\left\{\left| \frac{k_{i_0}^2(\lambda,t)d_{i_0}}{\epsilon(1+k_{i_0}(\lambda,t)d_{i_0})}\right| e^{-c\frac{d_{i_0}(x,t)}{\epsilon}},e^{-c\epsilon^{-\rho}}\right\}.\label{bdout0}\end{align}

In particular 
$$\left|Pv_\epsilon+\frac{\nabla_u W\left( v_\epsilon\right)}{\epsilon^2}\right|\leq C \hbox{ for some constant }C,$$

$$\left|Pv_\epsilon+\frac{\nabla_u W\left( v_\epsilon\right)}{\epsilon^2}\right|\to 0 \hbox{ as } \epsilon \to 0\hbox{ point-wise}$$
and
\begin{align}\int_0^t \int_{\Omega\setminus B_{\tilde{\delta}}}
 &\calH_{B_1}(x,y,t-s) \left|Pv_\epsilon+\frac{\nabla_u W\left( v_\epsilon\right)}{\epsilon^2}\right|(y,s)
dy ds \notag \\ = & \int_{t-\delta}^t \int_{\Omega\setminus B_{\tilde{\delta}}}
\calH_{B_1}(x,y,t-s) \left|Pv_\epsilon+\frac{\nabla_u W\left( v_\epsilon\right)}{\epsilon^2}\right|(y,s)
dy ds  \notag \\ &+\int_0^{t-\delta} \int_{\Omega\setminus B_{\tilde{\delta}}}
\calH_{B_1}(x,y,t-s) \left|Pv_\epsilon+\frac{\nabla_u W\left( v_\epsilon\right)}{\epsilon^2}\right|(y,s)
dy ds \notag \\
\leq & C\delta+ \frac{1}{\delta}\int_0^t \int_{\Omega\setminus B_{\tilde{\delta}}}
\left|Pv_\epsilon+\frac{\nabla_u W\left( v_\epsilon\right)}{\epsilon^2}\right|(y,s)
dy ds \to C\delta \hbox{ as }\epsilon \to 0. \label{bdout}\end{align}

\begin{rem} \label{unifbounds2}
As in Remark \ref{unifbounds} we have

$$ \int_{\Omega\setminus B_{\tilde{\delta}}}
\left|Pv_\epsilon+\frac{\nabla_u W\left( v_\epsilon\right)}{\epsilon^2}\right|(y,s)
\leq  C(e^{-c\epsilon^{\rho-1}}+\epsilon^{2\rho}) $$
\end{rem}

\begin{rem}\label{remss2}
As in Remark \ref{remss} when  $\sup_i|k_i|\leq f(t)$ where $f(t)$ is an integrable function of $t$ the computation above carries over in a similar way. 
Moreover,  the self-similar solutions in \cite{selsimos} satisfy 

\begin{align*}\int_0^t \int_{ B_{\tilde{\delta}- \epsilon}(O(t))\setminus B_{\epsilon^\rho}(O(t))}
 & \calH_{B_1}(x,y,t-s)  \left|Pv_\epsilon+\frac{\nabla_u W\left( v_\epsilon\right)}{\epsilon^2}\right|(y,s)
dy ds \\
\leq &C\delta+\frac{C}{\delta}\sqrt{T}\epsilon^{2\rho},\end{align*}
where the constants are independent of $\epsilon$ and depend linearly  on uniform bounds of $\sqrt{t}k_i$.

\end{rem}

\medskip

\noindent{\bf Transition Regions:}

\medskip 

Now we need to find bounds in the transition regions. As before, the computations in   $B_{\epsilon^\rho}(O(t))\setminus  B_{\frac{\epsilon^\rho}{2}}(O(t))$  are analogous to the ones in $B_{\tilde{\delta}}(O(t))\setminus B_{\tilde{\delta}- \epsilon}(O(t))$ and  we only present the calculations in the first set in detail. The computations that follow are similar to the ones in \cite{stationary}. 
 \medskip
 
 \item 
{\em In $B_{\epsilon^\rho}(O(t))\setminus  B_{\frac{\epsilon^\rho}{2}}(O(t))$  } we have
$\eta_2\left(\frac{r(x,t)}{2\epsilon}+1-\frac{\tilde{\delta}}{2\epsilon}\right)=1 $ and $v_\epsilon=\tilde{v}_\epsilon$. 

Recall that  $\{\xi^{int}_i\}$ is  a partition of unity, therefore
$$u_*\left(\frac{R_{\theta(t)}(x-O(t))}{\epsilon} \right)=\sum_{i=1}^6 \xi^{int}_{i}(\theta-\theta(t))
\left(u_*\left(\frac{R_{\theta(t)}(x-O(t))}{\epsilon} \right)\right).$$
Combining this equation with the definition of $\tilde{\phi}$, equations \equ{convu*dec} and \equ{convu*dec2} it is easy to see for $j\leq n\leq 2$ that 
\be \frac{\partial^n}{\partial^j t \partial^{n-j}x}\left(u_*\left(\frac{R_{\theta(t)}(x-O(t))}{\epsilon} \right)-\tilde{\phi}_\epsilon(x,t)\right)\leq  \frac{C}{\epsilon^2\left(1+\left(\frac{\epsilon^\rho}{\epsilon}\right)^m\right)} \leq C\epsilon^{m(1-\rho)-2}\label{difufi}\ee
for every $|x-O(t)|>\epsilon^\rho$.

Notice that we can write 
\begin{align*}\tilde{v}_\epsilon(x,t)= \tilde{\phi}_\epsilon(x,t)&+\eta_1\left(\frac{x-O(t)}{\epsilon^\rho}\right)
\left( u_*\left(\frac{R_{\theta(t)}(x-O(t))}{\epsilon} \right)-\tilde{\phi}_\epsilon(x,t)\right)
%\sum_{i=1}^3\left( \xi^{int}_{2i}(\theta-\theta(t))
%\left(u_*\left(\frac{R_{\theta(t)}(x-O(t))}{\epsilon} \right)\right.\right.
%\\ &\left.\left.-\zeta_{ii+1}\left(
%\frac{d_i(x,t)}{\epsilon}\right)\right)
%\frac{}{}\xi^{int}_{2i-1}(\theta-\theta(t))\left(u_*\left(\frac{R_{\theta(t)}(x-O(t))}{\epsilon} \right)-c_i\right)\right).
\end{align*}

Therefore
\begin{align}Pv_\epsilon+\frac{\nabla_u W(v_\epsilon)}{\epsilon^2}=&\left(P \tilde{\phi}_\epsilon(x,t)
+\frac{\nabla_u W(\tilde{\phi}_\epsilon)}{\epsilon^2}\right)\\& +P\left[\eta_1\left(\frac{x-O(t)}{\epsilon^\rho}\right)
\left( u_*\left(\frac{R_{\theta(t)}(x-O(t))}{\epsilon} \right)-\tilde{\phi}_\epsilon(x,t)\right)\right]
%\right. \\ &
%\sum_{i=1}^3 \left(\xi^{int}_{2i}(\theta-\theta(t))
%\left(u_*\left(\frac{R_{\theta(t)}(x-O(t))}{\epsilon} \right)-
%\zeta_{ii+1}\left(
%\frac{d_i(x,t)}{\epsilon}\right)\right)
%\frac{}{}\right. \notag \\
%&\left.\left.+\xi^{int}_{2i-1}(\theta-\theta(t))\left(u_*\left(\frac{R_{\theta(t)}(x-O(t))}{\epsilon} \right)-c_i\right)\right)\right]\notag \\
\\ &+\frac{\nabla_u W(\tilde{v}_\epsilon)}{\epsilon^2}-\frac{\nabla_u W(\phi_\epsilon)}{\epsilon^2} .\label{ect1}   \end{align}

Remark \ref{difbal} implies that we can bound $P \tilde{\phi}_\epsilon
+\frac{\nabla_u W(\tilde{\phi}_\epsilon)}{\epsilon^2}$   and its convolution with the heat kernel as before,  yielding
\begin{align}\int_0^t &\int_{ B_{\epsilon^\rho}(O(t)) \setminus B_{\frac{\epsilon^\rho}{2}(O(t))}}
 \calH_{B_1}(x,y,t-s) \left|P\tilde{\phi}_\epsilon+\frac{\nabla_u W\left( \tilde{\phi}\right)}{\epsilon^2}\right|(y,s)
dy ds \notag \\  
\leq & C\delta+ \frac{1}{\delta}\int_0^t \int_{ B_{\epsilon^\rho}(O(t)) \setminus B_{\frac{\epsilon^\rho}{2}(O(t))}}
\left|P\tilde{\phi}_\epsilon+\frac{\nabla_u W\left( \tilde{\phi}_\epsilon\right)}{\epsilon^2}\right|(y,s)
dy ds \to C\delta \hbox{ as }\epsilon \to 0. \label{bdt1}\end{align}

%$u_*(x)$ and its derivatives converge exponentially to $\tilde{\phi}$ as $|x| 
Noticing that
\begin{align*}\left|\frac{\nabla_u W(\tilde{v}_\epsilon)}{\epsilon^2}-\frac{\nabla_u W(\tilde{\phi}_\epsilon)}{\epsilon^2}\right|
\leq & C
\left|\frac{\tilde{v}_\epsilon-\tilde{\phi}_\epsilon}{\epsilon^2}\right|
\\ \leq & C \eta_1\left(\frac{x-O(t)}{\epsilon^\rho}\right)\frac{
\left|u_*\left(\frac{R_{\theta(t)}(x-O(t))}{\epsilon} \right)
-\tilde{\phi}_\epsilon(x,t)\right|}{\epsilon^2}.\end{align*}
Using equation \equ{difufi}
we have that the second and third terms of \equ{ect1} are bounded by
$C\epsilon^{\tilde{m}}$ for $\tilde{m}=m(1-\rho)-2$ and $m$ arbitrarily large, then we have that
$$\left|Pv_\epsilon+\frac{\nabla_u W(v_\epsilon)}{\epsilon^2}\right|\leq\left|P \tilde{\phi}
+\frac{\nabla_u W(\phi_\epsilon)}{\epsilon^2}\right|+ C\epsilon^{\tilde{m}},$$
for $m$ arbitrarily large.
This implies
\begin{align}\int_0^t & \int_{ B_{\epsilon^\rho}(O(t)) \setminus B_{\frac{\epsilon^\rho}{2}(O(t))}}
 \calH_{B_1}(x,y,t-s) \left|Pv_\epsilon+\frac{\nabla_u W\left(v_\epsilon \right)}{\epsilon^2}\right|(y,s)
dy ds \notag \\  
\leq & C\left(\delta+\frac{1}{\delta}\int_0^t \int_{ B_{\epsilon^\rho}(O(t)) \setminus B_{\frac{\epsilon^\rho}{2}(O(t))}}
\left|P\tilde{\phi}_\epsilon+\frac{\nabla_u W\left( \tilde{\phi}_\epsilon\right)}{\epsilon^2}\right|(y,s)+\epsilon^{\tilde{m}}\right) . \label{bdt2}\end{align}

\item {\em In $B_{\tilde{\delta}}(O(t))\setminus B_{\tilde{\delta}- \epsilon}(O(t))$} estimates are similar. Notice that for $x\in B_{\tilde{\delta}}(O(t))\setminus B_{\tilde{\delta}- \epsilon}(O(t))$ we have
$$v_\epsilon(x,t)= \left(1-\eta_2\left(\frac{r(x,t)}{2\epsilon}+1-\frac{\tilde{\delta}}{2\epsilon}\right)\right)\phi_\epsilon(x,t)
+\eta_2\left(\frac{r(x,t)}{2\epsilon}+1-\frac{\tilde{\delta}}{2\epsilon}\right)\tilde{\phi_\epsilon}(x,t).$$ Therefore

\begin{align*} Pv_\epsilon+\frac{\nabla_u W(v_\epsilon)}{\epsilon^2}=&\left(P \phi
+\frac{\nabla_u W(\phi_\epsilon)}{\epsilon^2}\right) \\
&+P\left(\eta_2\left(\frac{r(x,t)}{2\epsilon}+1-\frac{\tilde{\delta}}{2\epsilon}\right)\left(\tilde{\phi}_\epsilon-\phi_\epsilon \right)\right)\\
& +\left(\frac{\nabla_u W(v_\epsilon)}{\epsilon^2}-\frac{\nabla_u W(\phi_\epsilon)}{\epsilon^2}\right).
\end{align*}

 $P \phi_\epsilon
+\frac{\nabla_u W(\phi_\epsilon)}{\epsilon^2}$ can be bounded as before (see \equ{bdout}). Notice that by definition of $\phi_\epsilon$ and $\tilde{
\phi}_\epsilon$ we have
$\phi_\epsilon-\tilde{\phi}_\epsilon$ is equal to 0 near the interfaces. In fact, since $|x-O(t)|\geq \tilde{\delta}-\epsilon$, we have $(\phi-\tilde{\phi})(x,t)=0$ for every $x$ such that $d_i(x,t)\leq \min\{(\tilde{\delta} -\epsilon)\cos \delta_{int}, \delta \}$. Outside from these regions, both functions converge exponentially to the corresponding constant $c_i$ (see \equ{ecint3}), which implies
\begin{align*} \left|Pv_\epsilon+\frac{\nabla_u W(v_\epsilon)}{\epsilon^2}\right|\leq&\left|P \phi_\epsilon
+\frac{\nabla_u W(\phi_\epsilon)}{\epsilon^2}\right| 
+
 C\frac{e^{-c\epsilon^{\rho-1}}}{\epsilon^2} \end{align*}
and

\begin{align}\int_0^t \int_{B_{\tilde{\delta}}(O(t))\setminus B_{\tilde{\delta}- \epsilon}}  
 &\calH_{B_1}(x,y,t-s) \left|Pv_\epsilon+\frac{\nabla_u W\left(v_\epsilon \right)}{\epsilon^2}\right|(y,s)
dy ds \notag \\  
\leq & C\left(\delta+\frac{1}{\delta}\int_0^t \int_{B_{\tilde{\delta}}(O(t))\setminus B_{\tilde{\delta}- \epsilon}}
\left|P\phi_\epsilon(x,t)+\frac{\nabla_u W\left( \phi_\epsilon\right)}{\epsilon^2}\right|(y,s)+C\frac{e^{-c\epsilon^{\rho-1}}}{\epsilon^2} \right)\\
&\to C\delta \hbox{ as } \epsilon \to 0 . \label{bdt3}\end{align}

\medskip

Finally, combining \equ{bd1}, \equ{bd2}, \equ{bdout} , \equ{bdt2} and  \equ{bdt3}   we conclude 
\begin{align}\int_0^t \int_{\Omega}  
 &\calH_{B_1}(x,y,t-s) \left|Pv_\epsilon+\frac{\nabla_u W\left(v_\epsilon \right)}{\epsilon^2}\right|(y,s)
dy ds 
\to C\delta \hbox{ as }\epsilon \to 0.\end{align}
Since $\delta$ is arbitrary this concludes the proof.

\end{itemize}
\end{proof}

\begin{rem}\label{cotfinal}
It is easy to see from the proof above and Remarks \ref{unifbounds} and \ref{unifbounds2}
that there is an $l>0$ such that
$$|F_\epsilon(w_\epsilon,\psi_\epsilon)-w_\epsilon|\leq C\left( \delta+
\frac{\epsilon^{2l}}{\delta}\right),$$
for every $\delta>0$. By choosing $\delta=\epsilon^l$ we have
$$|F_\epsilon(w_\epsilon,\psi_\epsilon)-w_\epsilon|\leq C\epsilon^{l}.$$
This inequality combined with Corollary \ref{corcotprinc}
implies
\be \sup_{\Omega\times[0,T]}|u_\epsilon-v_\epsilon|\leq C\epsilon^l. \ee

From remarks \ref{remss} and \ref{remss2} we have that this bound also holds for the self-similar solutions in \cite{selsimos}.
\end{rem}

\section{Proof of Corollaries \ref{cor1} and \ref{cor2}}\label{pfcor}

\medskip
\noindent{\bf Proof of Corollary \ref{cor1}}

Let  $\calT_1(\cdot,t)=\{\gamma_1^i(\cdot)\}_{i=1}^3$ and $\calT_2(\cdot,t)=\{\gamma_2^i(\cdot)\}_{i=1}^3$ 
Suppose that we have $\calT_1(\cdot,0)=\{\sigma_1^i(\cdot)\}_{i=1}^3=\{\sigma_2^i(\cdot)\}_{i=1}^3=\calT_2(\cdot,0)$. Moreover, we assume $\sigma_1^i(\cdot)=\sigma_2^i(\cdot)$ for every $i$.
Let $d_j^i(x,t) $ be the distance to $\gamma_j^i(\cdot,t)$. 
Consider $v_\epsilon^1(x,t)$ and $v_\epsilon^2(x,t)$ to be the functions defined by Theorem \ref{mainteo}
for the triods $\calT_1$ and $\calT_2$ respectively and  let $u_\epsilon^1(x,t)$ and $u_\epsilon^2(x,t)$
be the associated solutions.

Notice that since the distance functions are independent of the parametrization of $\sigma_j^i$. Hence, we have that $d_1^i(x,0)=d_2^i(x,0)$ and
 the definitions of $v_\epsilon^j$ imply that
$\psi_\epsilon^1(x)=v_\epsilon^1(x,0)=v_\epsilon^2(x,0)=\psi_\epsilon^2(x)$.  Uniqueness gives us
$$u_\epsilon^1(x,t)=u_\epsilon^2(x,t),$$ for every $t>0$ and $\epsilon>0$.
This implies that for every $t>0$ 
$$\{ x \in \Omega: \lim_{\epsilon \to 0} u_\epsilon(x,t)^1\ne c_k \}= \{ x \in \Omega: \lim_{\epsilon \to 0} u_\epsilon^2(x,t)\ne c_k \}.$$

Or equivalently $\calT_1(\cdot,t)=\calT_2(\cdot,t)$ for every $t>0$.$\Box$

\medskip

\noindent {\bf Proof of Corollary \ref{cor2}}

Consider the triod $\calT_n(x,t)=\left\{\frac{1}{\beta_n}\gamma^i(\cdot,\beta_n^2 t)\right\}$ (as defined in the statement of Corollary \ref{cor2}).  It is easy to verify that 
$\calT_n$ satisfies equation \equ{triod}.  Consider any symmetric potential $W$ satisfying the conditions of Theorem \ref{mainteo}.
Let $u_\epsilon^n$ and $v_\epsilon^n$ be the solution and the approximation defined by Theorem \ref{mainteo} for $\calT_n$. Similarly, let  $\calT_{self}$ be the self-similar triod defined in \cite{selsimos} with initial condition  $\left\{ \lambda \frac{\gamma_\lambda^i(0,t)}{|\gamma_\lambda^i(0,t)|}, \lambda \in \rr \right\}_{i=1}^3$ and let  $u_\epsilon^{self}$, $v_\epsilon^{self}$ be respectively the solution and the approximation defined by Theorem \ref{mainteo} for $\calT_{self}$. 

Notice that the curvatures $k_i^n$ of the curves of the triod $\calT_n$ satisfy
 $$\sup_i |k_i^n|(x,t)=\sup_i \beta_n  |k_i|(x, \beta_n^2 t)\leq \beta_n \frac{C}{\beta_n \sqrt{t}}=\frac{C}{ \sqrt{t}}.
$$
Fix a ball $B_R$ such that $O(t)\in B_R$ for $t\in [0,T]$. Theorem \ref{mainteo} implies
$$\sup_{B_{2R}\times [0,T]}|u_\epsilon^n-v_\epsilon^n|\to 0 \hbox{ as  }\epsilon \to 0$$
and
 $$\sup_{B_{2R}\times [0,T]}|u_\epsilon^{self}-v_\epsilon^{self}|\to 0 \hbox{ as  }\epsilon \to 0.$$
 Moreover,  Remark \ref{cotfinal} implies that
    $$\sup_{B_{2R}\times [0,T]}|u_\epsilon^n-v_\epsilon^n|\leq C\epsilon^{l}$$
and
 $$\sup_{B_{2R}\times [0,T]}|u_\epsilon^{self}-v_\epsilon^{self}|\leq C\epsilon^{l},$$
where $C$ and $l$ are independent of $n$.

Consider now a smooth function $\chi$ supported in $B_{2R}$ such that $\chi\equiv 1$ in $B_R$.
It is easy to see that for every $\epsilon>0$  the functions $\chi(x)u_\epsilon^n(x,t)$ and $\chi(x)u_\epsilon^{self}(x,t)$ satisfy the same nonlinear parabolic equation
in $B_{2R}$ with boundary values equal to 0 for every $t>0$. 

Notice that since $\gamma^i(0,0)=0$ we have that $\calT_n(0,t)\to  \left\{ \lambda \frac{\gamma_\lambda^i(0,t)}{|\gamma_\lambda^i(0,t)|}, \lambda \in \rr_+\right\}_{i=1}^3$ uniformly in compact subsets.
By the definition of the initial condition of Theorem \ref{mainteo} , holds
$$\chi(x)u_\epsilon^n(x,0)\to \chi(x)u_\epsilon^{self}(x,0) \hbox{ as }n\to \infty$$ uniformly in $B_{2R}$.
Standard parabolic theory implies for every $\epsilon>0$ that
$$\chi(x)u_\epsilon^n(x,t)\to \chi(x)u_\epsilon^{self}(x,t) \hbox{ as }n\to \infty.$$
uniformly in $B_{2R}\times [0T]$
We conclude that
\begin{align*}\sup_{B_R\times [0,T]}|v_\epsilon^n-v_\epsilon^{self}|(x,t)=&\sup_{B_R\times [0,T]}|\chi(v_\epsilon^n-v_\epsilon^{self})|(x,t))\\
\leq & \sup_{B_{2R}\times [0,T]}|u_\epsilon^n-v_\epsilon^n|(x,t) +\sup_{B_{2R}\times [0,T]}|\chi(u_\epsilon^n-u_\epsilon^{self})|(x,t)\\ &+  \sup_{B_{2R}\times [0,T]}|u_\epsilon^{self}-v_\epsilon^{self}| (x,t)\\
\leq& 2C\epsilon^{2l}+sup_{B_{2R}\times [0,T]}|\chi(u_\epsilon^n-u_\epsilon^{self})|(x,t).
\end{align*}

Taking $n\to \infty$ and $\epsilon \to 0$ we conclude that
$$\lim_{\epsilon\to 0} \lim_{n \to \infty}\sup_{B_R\times [0,T]}|v_\epsilon^n-v_\epsilon^{self}|(x,t)=0.$$
The definition of  $v_\epsilon^n$ implies the result. $\Box$

\section{Further Comments}\label{finsec}

We would like to finish this paper with some open problems that we hope in the future 
can be answered using this representation.
\medskip
\begin{enumerate}
\item{\bf Long time existence:}
In the work by Mategazza, Novaga and Tortorelli \cite{moybycm} singularities type I were excluded. However this was not possible for singularities type II. From the point of view of the Allen-Cahn equation it is easy to see that for $\epsilon>0$ solutions to \equ{laeq}-\equ{ci}-\equ{bc} exist for all times. Moreover, standard computations show that as $t\to \infty$ this solutions converge to solution of the associated elliptic equation. Work of Baldo \cite{minintsb}
shows that solutions to the stationary problem sub-converge to piecewise constant function with interfaces that are piecewise linear. This exactly corresponds to the long time behavior expected for solutions to \equ{triod}. It would be interesting to use the relaxation method presented in this paper to prove long time existence as described above.

\medskip

\item{\bf Convergence to Self-Solutions in the general case:}
Corollary \ref{cor2} showed that any solution to \equ{triod} that satisfies 
$|k_i|\leq \frac{C}{\sqrt{t}}$ converges to a self-similar solution described by \cite{selsimos}. We expect that this bound on the curvatures can be removed.

\medskip

\item{\bf Existence of solution with arbitrary initial data:}
 Since solutions to  \equ{laeq}-\equ{ci}-\equ{bc} exist for $\epsilon>0$ for arbitrary initial conditions (in particular relaxation of triods satisfying arbitrary angle conditions), this representation might provide a weak solution in these cases (by observing the interfaces formed as $\epsilon\to 0$).
  It would interesting to understand  such weak solutions.

\end{enumerate}

\bibliography{mariel} 
\end{document}